\documentstyle[12pt]{article}

\newtheorem{theorem}{Theorem}[section]
\newtheorem{lemma}[theorem]{Lemma}
\newtheorem{corollary}[theorem]{Corollary}

\newcommand{\R}{{\bf R}}

\newcommand{\A}{{\cal A}}
\newcommand{\F}{{\cal F}}
\newcommand{\ch}{\raisebox{.4ex}{$\chi$}}
\renewcommand{\sp}{\vspace{1ex}}

\newcommand{\pl}{\partial}
\newcommand{\nm}{\parallel}
\newcommand{\ov}{\over}
\newcommand{\ra}{\rightarrow}
\newcommand{\ph}{\varphi}
\newcommand{\qed}{\hfill $\Box$}
\newcommand{\proof}{{\em Proof. }}
\newcommand{\be}{\begin{equation}}
\newcommand{\ee}{\end{equation}}
\newcommand{\bq}{\begin{equation}}
\newcommand{\eq}{\end{equation}}
\newcommand{\nn}{\nonumber}
\newcommand{\ba}{\begin{array}}
\newcommand{\ea}{\end{array}}

\newcommand{\inv}{^{-1}}
\newcommand{\iy}{\infty}
\newcommand{\al}{\alpha}

\newcommand{\la}{\lambda}

\newcommand{\ind}{\mbox{\rm ind\,}}
\newcommand{\tr}{\mbox{\rm tr\,}}

\newcommand{\tf}{\tilde f}

\begin{document}

\title{\bf Determinants of Airy Operators and Applications to Random
Matrices}
\author{Estelle L. Basor\thanks{Supported in part by NSF Grant
DMS-9623278.}\\
               \emph {Department of Mathematics}\\
               \emph {California Polytechnic State University}\\
               \emph {San Luis Obispo, CA 93407, USA}
        \and
      Harold Widom\thanks{Supported in part by NSF Grant DMS-9732687.}\\
               \emph {Department of Mathematics}\\
                   \emph {University of California}\\
                   \emph {Santa Cruz, CA 95064, USA}}
\date{}

\maketitle

\begin{abstract}
The purpose of this paper is to describe asymptotic formulas for
determinants of certain operators that are analogues of
Wiener-Hopf operators. The determinant formulas yield information
about the distribution functions for certain random variables that arise
in
random matrix theory when one rescales  at ``the edge of the
spectrum''.
\end{abstract}


\section{Introduction}

This paper is concerned with the asymptotics of Fredholm determinants of
operators that arise naturally in random matrix theory and
are similar in many ways to finite Wiener-Hopf operators.
The operators, denoted by $A_{\alpha}(f),$ are
integral operators on $L^{2}(\bf{R})$ with kernel given by
\bq
f(x/\alpha)\int_{0}^{\iy}A(x+z)A(z+y)dz\label{int.0}
\eq
where
\[
A(x)  =  \frac{1}{2\pi}\int_{-\iy}^{\iy}e^{it^{3}/3}e^{itx}dt,
\]
and $f \in L^{\iy}(\bf{R}).$ The function $A(x)$ is the Airy function,
generally denoted
by Ai$(x)$, and for this reason we call our operators $A_{\al}(f)$ {\em
Airy operators}.
We will refer to the function $f$ as the
{\it symbol} of the Airy operator.

If the term
$\int_{0}^{\iy}A(x+z)A(y+z)dz$
in (\ref{int.0}) is replaced by the sine kernel
\[ \frac{\sin\pi(x-y)}{x-y}, \] the resulting operator has the same
Fredholm
determinant as a finite
Wiener-Hopf operator, whose asymptotics are very well known \cite{K}.
The similarities with the Wiener-Hopf (i.e., sine kernel) case become
less
surprising after the observation that
\bq
\int_{0}^{\iy}A(x+z)A(y+z)dz  =  \frac{A(x)A'(y) - A(y)A'(x)}{x-y}. \nn
\eq
The proof of this as well as many other facts about the above kernel
can be found in \cite{TW}. We shall use the fact that $A(x)$
 is rapidly decreasing at $+\iy$ and
$O(|x|^{-1/4})$ at $-\iy$. For the complete asymptotics
of the Airy function, we refer the reader to \cite{E}.

Under appropriate conditions $A_{\alpha}(f)$ is a trace class operator,
and thus the
Fredholm determinant
$$\det(I + A_{\alpha}(f))$$
is defined. The main goal of the paper is to compute the asymptotics
of this determinant as $\alpha \rightarrow \iy.$

The motivation for finding an asymptotic formula for the Fredholm
determinant comes from random matrix theory, in studying so-called
{\em linear statistics}, which are certain functions of
the eigenvalues of random matrices. After a rescaling
at ``the edge of the spectrum'' their characteristic functions become
Fredholm determinants of our Airy operators. For general information
about
random matrices, we refer the reader to \cite{M}. For information
about the connection of random matrices, characteristic functions
and the Airy operators we refer the reader to \cite{B} and \cite{TW}.

The paper is organized as follows. In the second section we derive the
basic
properties of $A_{\alpha}(f)$ and related operators.  In the third
section we prove
through a series
of lemmas that for appropriate functions $f$ and $F$
\[\lim_{\al\ra\iy}\tr[F(A_{\al}(f)) -A_{\al}(F\circ f)]=\tr[F(W(g)) -
W(F \circ g)],\]
where $W(g)$ is the Wiener-Hopf operator with symbol $g(x)=f(-x^2)$.
(The precise definition of $W(g)$ will be given at the end of the next
section.)
The trace of the second
operator on the left is easy to compute asymptotically. Taking $F(z) =
\log(1+z)$
and using the known formula for the trace
on the right side, we find that the Fredholm determinant is given
asymptotically as
$\al \ra \iy$ by
\bq
\det\,(I + A_{\alpha}(f)) = \exp \left\{ c_{1}\,\alpha^{3/2} +
c_{2}+o(1)\right\},\label{int.2}
\eq
 where
\[
c_{1}  =  {1\over\pi}\,\int_{0}^{\iy}\sqrt{x}\,\log(1+f(-x))\,dx,
\]
\[
c_{2}  =  \frac{1}{2}\int_{0}^{\iy}x\,(G(x))^{2}\,dx,
\]
and $$G(x) = \frac{1}{{2\pi}}\int_{-\iy}^{\iy}e^{ixy}\,\log\,(1
+f(-y^{2}))\,dy.$$
This is proved under the assumption that $f$ is a Schwartz function
(although we could get by with much less) and $1+f(x)\ne0$ for $x\le0$.

This formula bears a strong resemblance to the corresponding asymptotic
formula in the
classical finite Wiener-Hopf case. The most notable difference is the
power
$\al^{3/2}$ in the first term of the asymptotics.

In the last section we describe the implications of
formula (\ref{int.2}) for random matrix theory. The formula, as in the
analogous Wiener-Hopf or Bessel kernel case (see \cite {B} for details),
proves that the distribution functions for certain linear statistics,
now scaled
at the edge of the spectrum, are asymptotically Gaussian. The
recurrence of the Gaussian distribution highlights the
universality seen again and again in random matrix models.

\section{Basic properties of the Airy operator}

We begin by defining the {\it Airy transform} $\A$. For $g\in L^{2}(\R)$
we define $\A(g)$  by the formula
$$\A (g) = {\cal F}^{-1}{M}_{h} {\cal F}^{-1}(g),$$
where ${\cal F}(g)(x)$ is the Fourier transform of $g$ given by
$$\frac{1}{\sqrt{2\pi}}\int_{-\iy}^{\iy}g(t)e^{-ixt}dt,$$ ${\cal
F}^{-1}$ is
the inverse transform and ${M}_{h}$ denotes multiplication by the
function
$h(t) = e^{it^{3}/3}.$ We
will also use the standard notation $\hat{g}$ and $\check{g}$ for the
Fourier transform and inverse transforms respectively. Observe
that for $g\in L^{1}\cap L^{2}$ we have
$\A (g)(x) = \int_{-\iy}^{\iy}A(x+y)g(y)dy$ and,
just as in the Fourier transform case, $\A(g)$ is the $L^{2}$  limit
of $\int_{-B}^{B}A(x+y)g(y)dy$ as $B \ra \iy$ for all $g\in L^2$.

\begin{lemma}\label{2.1}
The Airy transform is unitary on $L^{2}$ and satisfies $\A \inv = \A.$
\end{lemma}
\proof Clearly $\A$ is unitary since ${\cal F}$ and $M_h$ are.
The Fourier inversion formula says that ${\cal F}^{-1}={\cal J}\,{\cal
F}=
{\cal F}\,{\cal J}$, where ${\cal J}g(x)=g(-x)$. It follows that
$$\A={\cal F}\,{\cal J}\,M_h\,{\cal J}\,{\cal F}={\cal
F}\,M_{h^{-1}}\,{\cal F}=\A^{-1}.$$

\vspace{-4.3ex}\qed

\sp

Given the above definition of the Airy transform, we see
that the Airy operator defined by (\ref{int.0}) is alternatively defined
as the operator $M_{f_{\alpha}} \A P \A,$ where $P$ is multiplication by
$\ch_{\R^{+}}$ and $f_{\alpha}(x) = f(x/\alpha).$ It is this
representation of the operator which we will use.
For appropriate $f$ this operator in turn will have the same Fredholm
determinant
as $\A M_{f_{\alpha}} \A P.$
We next derive a representation for the kernel of $\A M_{f_{\alpha}} \A$
for
a class of functions $f$.

\begin{lemma}\label{2.2} Suppose that $f$ is the inverse Fourier
transform of a finite
measure $\mu$,
$$ f(x) = {1\over\sqrt{2\pi}}\int_{-\iy}^{\iy}e^{i\xi x}\,d\mu(\xi).$$
Then the kernel of the operator $\A M_{f} \A$ is given by
the formula
\bq
\frac{1}{\sqrt{8}\pi}\int_{-\iy}^{\iy}\frac{e^{-i\xi^{3}/12}}{\sqrt{i\xi
+ 0}}\,
e^{-i(x+y)\xi/2}\,
\,e^{i(x-y)^{2}/4\xi}\,d\mu(\xi), \label{pro.1}
\eq
where $\sqrt{i\xi+0}$ is defined by taking $\arg(i\xi+0) $ equal
to $\pi/2$ when $\xi > 0$ and $-\pi/2$ when $\xi < 0$.
\end{lemma}
\proof
Consider first the case where $\mu$ is a unit mass at the point $\eta$,
and let $K_{\eta}$
be the operator with the corresponding kernel (\ref{pro.1}). Then $f(x)$
is the
function
$e_{\eta}(x)=e^{i\eta x}/\sqrt{2\pi}.$
By Lemma 2.1 we see that we have to show
\[K_{\eta}=\F^{-1}\,M_h\,\F^{-1}\,M_{e_{\eta}}\,\F\,M_{h^{-1}}\,\F,\]
or equivalently
\[K_{\eta}\,\F^{-1}=\F^{-1}\,M_h\,\F^{-1}\,M_{e_{\eta}}\,\F\,M_{h^{-1}}.\]
Let us compute both sides applied to a function in $L^2$.

Notice first that $\F^{-1}\,M_{e_{\eta}}\,\F$ takes a function
$\ph(\xi)$ into
$\ph(\xi+\eta)/\sqrt{2\pi}$. Therefore \linebreak
$M_h\,\F^{-1}\,M_{e_{\eta}}\,\F\,M_{h^{-1}}$
takes
$\ph(\xi)$ into
\[{1\over\sqrt{2\pi}}\,e^{i\xi^3/3}\,e^{-i(\xi+\eta)^3/3}\,\ph(\xi+\eta)
={1\over\sqrt{2\pi}}\,e^{-i(\xi^2\eta+\xi\eta^2+\eta^3/3)}\,\ph(\xi+\eta).\]
Hence
\[\F^{-1}\,M_h\,\F^{-1}\,M_{e_{\eta}}\,\F\,M_{h^{-1}}\,\ph(x)
={1\over2\pi}\,\int_{-\iy}^{\iy}e^{i\xi
x}\,e^{-i(\xi^2\eta+\xi\eta^2+\eta^3/3)}\,
\ph(\xi+\eta)\,d\xi\]
\bq ={e^{-i(\eta x+\eta^3/3)}\over 2\pi}\,
\int_{-\iy}^{\iy}e^{i\xi
x}\,e^{-i(\xi^2\eta-\xi\eta^2)}\,\ph(\xi)\,d\xi.\label{A}\eq

On the other hand, we have
\[K_{\eta}\,\F^{-1}\,\ph(\xi)=\frac{1}{\sqrt{8}\pi}
\frac{e^{-i\eta^{3}/12}}{\sqrt{i\eta + 0}}\,\frac{1}{\sqrt{2\pi}}\,
\int_{-\iy}^{\iy}\,\ph(\xi)\,d\xi\,
\int_{-\iy}^{\iy}\,e^{i\xi
y}\,e^{-i(x+y)\eta/2}\,e^{i(x-y)^{2}/4\eta}\,dy.\]
(If $\ph$ is a Schwartz function, for example, the interchange of order
of integration
this involves can be justified by integration by parts, and it suffices
to show our
two operators agree when applied to Schwartz functions.) The inner
integral
is easily computed and found to equal
\[2\sqrt{\pi}\,\sqrt{i\eta+0}\,e^{i(\xi-\eta)x}\,e^{-i(\xi^2\eta-\xi\eta^2)}\,
e^{-i\eta^{3}/4}.\]
Thus we see that $K_{\eta}\,\F^{-1}\,\ph(\xi)$ is equal to the right
side of (\ref{A}).

The lemma is established for the special case of a unit point mass, and
so for any linear
combination of these. To establish the general result we approximate our
given $\mu$ by
a sequence $\{\mu_n\}$ of linear combinations of point masses such that
the
measures $\mu_n$ are uniformly bounded and
$\int \ph(\xi)\,d\mu_n(\xi)\ra \int \ph(\xi)\,d\mu(\xi)$ for any
function $\ph$ which
is bounded and continuous. Then the corresponding functions $f_n$
converge boundedly and pointwise to $f$ and so the operators $\A M_{f_n}
\A$
converge strongly to $\A M_{f} \A$. For the corresponding operators
$K_n$, it is easy to
see that for Schwartz functions $g_1$ and $g_2$ we have
$(K_n\,g_1,\,g_2)\ra
(K\,g_1,\,g_2)$, so $K_n\ra K$ weakly. Hence, since $\A M_{f_n} \A=K_n$
for
each $n$, we have $\A M_{f} \A=K$.
\qed\sp

To end this section we recall the definition  of a Wiener-Hopf operator
and certain of
its properties. For a function $g \in L^{\iy}(\R)$ the operator $W(g)$
on
$L^{2}(\R^{+})$ (which we identify with the functions in $L^{2}(\R)$
which vanish on $\R^-$) is defined by
$$W(g)= P \F^{-1} M_{g} \F P.$$
This is the {\em Wiener-Hopf operator with symbol $g$}. Notice the
analogy with the operators
$P\A M_{f}\A P$. The fact that there is more than just an analogy will
become apparent in
the next section. One often sees a Wiener-Hopf operator defined as an
operator on
$L^{2}(\R^{+})$ with kernel of the form $k(x-y)$ where $k\in L^1(\R)$.
This operator is equal to $W(g)$ with $g(x)=\int e^{-ixu}\,k(u)\,du$.\sp

We state as a lemma two basic facts about Wiener-Hopf operators.
\begin{lemma} \begin{itemize}
\item [a)] The spectrum of $W(g)$ is contained in the convex hull of the
essential
range of $g$.
\item [b)] If $g$ is continuous and $g(\pm\iy)=0$ then $\la$ is not
in the spectrum
of $W(g)$ if and only if $\la\ne0,\ g-\la\ne0$ and
\[i(g-\la):={1\ov2\pi}\raisebox{-1.5ex}{$\stackrel{\textstyle
\Delta}{\scriptstyle -\iy<x<\iy}$}
 \arg\,(g(x)-\la)=0.\]
 \end{itemize}
 \end{lemma}

\section{Trace norm estimates and the Airy limit theorem}

We assume from now on that $f$ is a Schwartz function.
The reader can verify that this requirement is
too restrictive and can, for example, be replaced by a weighted space
condition. However,
assuming that $f$ is a Schwartz function
simplifies the proofs and increases the clarity of the arguments.

Recall that the Airy operator $A_{\al}(f)$ equals $M_{f_{\alpha}}\A P
\A$ and is
thus similar to $\A M_{f_{\alpha}} \A P$, which in turn is to unitarily
equivalent to the operator $U^{-1}\A M_{f_{\alpha}} \A P U$
where $U$ is the unitary operator defined by $Ug(x) =
\alpha^{1/4}g(x\sqrt{\alpha}).$ Note that $U$ commutes with $P$.
This operator will act as a replacement for the Airy operator
in the final computations. The next lemma, which involves a
modification of the above operator, will be important in
those computations.
\begin{lemma}\label{3.1}
The operator $U^{-1}(I-P) \A M_{f_{\alpha}} \A P U$ converges in the
trace norm
to the operator with kernel
$$\frac{1}{\sqrt{8}\pi}\int_{-\iy}^{\iy}\frac{\hat{f}(\xi) }{\sqrt{i\xi
+ 0}}\,
e^{i(x-y)^{2}/4\xi}\,d\xi\,\ch_{\R^{-}}(x)\,\ch_{\R^{+}}(y)$$
as $\al \ra \iy$.
\end{lemma}
\proof
By Lemma 2.2 the kernel of the operator
$ (I-P) \A M_{f_{\alpha}} \A P $ is given by the formula
$$\frac{1}{\sqrt{8}\pi}\int_{-\iy}^{\iy}\frac{\al\,\hat{f}(\al
\xi)}{\sqrt{i\xi+0}}\,e^{-i\xi^3/12}\,e^{-i(x+y)\xi/2}\,
e^{i(x-y)^2/4\xi}\,d\xi\, \ch_{\R^{-}}(x)\ch_{\R^{+}}(y),$$
and thus the kernel of $U^{-1}(I- P) \A M_{f_{\alpha}} \A P U$ is given by
$$\frac{1}{\sqrt{8}\pi}\int_{-\iy}^{\iy}\frac{\hat{f}(
\xi)}{\sqrt{i\xi+0}}\,e^{-i\xi^3/12\al^{3}}\,e^{-i(x+y)\xi/2\al^{3/2}}\,
e^{i(x-y)^2/4\xi}\,d\xi\, \ch_{\R^{-}}(x)\ch_{\R^{+}}(y).$$ Changing
$x$ to $-x$ in the kernel for convenience gives
$$\frac{1}{\sqrt{8}\pi}\int_{-\iy}^{\iy}\frac{\hat{f}(
\xi)}{\sqrt{i\xi+0}}\,e^{-i\xi^3/12\al^{3}}\,e^{i(x-y)\xi/2\al^{3/2}}\,
e^{i(x+y)^2/4\xi}\,d\xi\, \ch_{\R^{+}}(x)\ch_{\R^{+}}(y).$$

We shall show that replacing by 1 each of the two factors in the
integrand which involve $\al$ leads to an error which is the kernel of
an
operator having trace norm $o(1)$.
We will use the general fact that the trace norm of an operator with
kernel $K(x,y)$, where
$y$ is confined to a set $J$, is at most $\nm K\nm_2+|J|\,\nm  \pl K/
\pl y\nm_2$, where the
norms are Hilbert-Schmidt norms.

We first look at the error kernel arising from the replacement
$e^{-i(x-y)\xi/2\al^{3/3}}\ra 1$, which is
$$   \frac{1}{\sqrt{8}\pi}\int_{-\iy}^{\iy}\frac{\hat{f}(
\xi)}{\sqrt{i\xi+0}}\,e^{-i\xi^3/12\al^{3}}\,(e^{i(x-y)\xi/2\al^{3/2}}-1)\,
e^{i(x+y)^2/4\xi}\,d\xi\, \ch_{\R^{+}}(x)\ch_{\R^{+}}(y).$$
This we call $K(x,y)$ and find
bounds on $K$ and $ \pl K/\pl y.$

Clearly
$K(x,y) = O(|x-y|/\alpha^{3/2})=O(w/\alpha^{3/2})$, where $w=x+y$.
We use this estimate for
$w\leq 1.$ To get a better estimate for $w\geq 1$ we
write the kernel as a constant times
$$\frac{1}{w^{2}}\int_{-\iy}^{\iy}\frac{\hat{f}(
\xi)}{\sqrt{i\xi+0}}\,e^{-i\xi^3/12\al^{3}}\,(e^{i(x-y)\xi/2\al^{3/2}}-1)\,
\xi^{2}\frac{d}{d\xi}e^{i(x+y)^2/4\xi}\,d\xi\,$$
and integrate by parts to obtain a constant times
$$\frac{1}{w^{2}}\int_{-\iy}^{\iy}\frac{d}{d\xi}\Big[\xi^{2}\frac{\hat{f}(
\xi)}{\sqrt{i\xi+0}}\,e^{-i\xi^3/12\al^{3}}\,(e^{i(x-y)\xi/2\al^{3/2}}-1)\Big]
\,
e^{i(x+y)^2/4\xi}\,d\xi.$$
Of course we apply the product rule. Differentiating the various factors
in the brackets
leads to an extra factor $\al^{-3}$ or $w\,\al^{-3/2}$, aside
from the external factor $1/w^2$.
If $w\geq 1$ we see therefore
that integration by
parts yields a factor $\al^{-3/2}<w>^{-1}=\al^{-3/2}(1+w^2)^{-1/2}$
assuming that
of course $f$ is in a Schwartz function. Integrating by parts once
more leads to  a factor of $\al^{-3}<w>^{-2}.$ Thus for all positive
$x$ and $y$ ($w\geq 1$ or $w\leq 1$)  we see that our kernel satisfies
$K(x,y)=O(\al^{-3/2}<w>^{-2}).$

We also have to estimate $\pl K(x,y)/\pl y.$
If we differentiate
$$\int_{-\iy}^{\iy}\frac{\hat{f}(
\xi)}{\sqrt{i\xi+0}}\,e^{-i\xi^3/12\al^{3}}\,(e^{i(x-y)\xi/2\al^{3/2}}-1)\,
e^{i(x+y)^2/4\xi}\,d\xi\,$$ we are left with two integrals. In one
integral
we get an
extra $\xi/\al^{3/2}$ in the integrand, and the factor
$e^{i(x-y)\xi/2\al^{3/2}}-1$ is replaced by $e^{i(x-y)\xi/2\al^{3/2}}$.
As
before this can be
seen to be  $O(\al^{-3/2}<w>^{-2}).$ In the other integral
we get an extra $w/\xi$ which changes the factor
$e^{i(x-y)\xi/2\al^{3/2}}-1$ to
$$w\frac{e^{i(x-y)\xi/2\al^{3/2}}-1}{\xi}.$$ Here after three
integration by parts we arrive at an estimate of
$O(\al^{-3/2}<w>^{-2}).$

We have shown that
\[K(x,y),\ \ {\pl\over\pl y}K(x,y)=O\left({\al^{-3/2}\over
1+x^2+y^2}\right).\]
If we use the general trace norm estimate stated above,
taking $J=(k,\,k+1)$
for $k=0,\; 1,\cdots$ and adding, we find that the trace
norm of the error operator is $O(\al^{-3/2})$.

If we consider the error due to the  replacement $e^{-i\xi^3/12\al^3}\ra
1$ the argument is essentially the same and we find a bound
for the trace norm of the resulting kernel of $O(\al^{-3}).$ This
completes the proof.
\qed\sp

Here and below we shall use the notations $O_1(\cdot)$ resp.
$o_1(\cdot)$ to denote families
of operators depending on the parameter $\al$ whose trace norms are
$O(\cdot)$ resp.
$o(\cdot)$.

\begin{lemma}. We have $P \A M_{f_{\alpha}} \A P=O_1(\al^{3/2})$ in
general and
$P \A M_{f_{\alpha}} \A P=o_1(1)$ if $f$ vanishes on $\R^-$.
\end{lemma}
\proof
The kernel of our operator on $L^2(\R^+)$ equals
$$\int_{-\iy}^{\iy}f(z/\al)\,A(x+z)\,A(y+z)\,dz.$$
For fixed $z$ the kernel $f(z/\al)\,A(x+z)\,A(y+z)$ is a
separable rank one kernel. To compute its trace norm, observe that
by the estimates on the Airy function we have
\[\int_0^{\iy}A(x+z)^2\,dx=\left\{
\begin{array}{ll}O(e^{-z}),&z>0,\\&\\O(<z>^{1/2}),&z<0.
\end{array}\right.\]
It follows that the trace norm of our operator is at most a constant
times
\[\int_{-\iy}^0|f(z/\al)|\,<z>^{1/2}\,dz+\int_0^{\iy}|f(z/\al)|\,e^{-z}\,dz,\]
and the assertions of the lemma follow easily.\qed

\begin{corollary}
$A_{\al}(f)$ is a trace class operator.
\end{corollary}
\proof
After the replacement $x\ra-x$ the kernel in the statement of Lemma~3.1
becomes
a Hankel operator with smooth kernel and thus is well-known to be trace
class. The lemma implies that
$U^{-1}(I-P) \A M_{f_{\alpha}} \A P U$
is trace class. Thus
$(I-P) \A M_{f_{\alpha}} \A P$ is trace class. Lemma~3.2 tells us that
$P\A M_{f_{\alpha}} \A P$ is trace class. Hence so is $\A M_{f_{\alpha}}
\A P$,
and $A_{\al}(f)=M_{f_{\alpha}} \A P \A$ is unitarily equivalent to this.
\qed\sp

We remark that this argument could have been made much earlier.
However it would have involved the same sort of estimates as in the
proof
of Lemma 3.1 and there was no point in doing this twice.\sp

The operator with kernel
$$\frac{1}{\sqrt{8}\pi}\int_{-\iy}^{\iy}\frac{\hat{f}(\xi) }{\sqrt{i\xi
+ 0}}
e^{i(x-y)^{2}/4\xi}d\xi$$ is a convolution operator, its
kernel is a function of $x-y$. We will denote it by $K_{f}.$
Thus $P K_{f} P$ is a Wiener-Hopf operator with symbol
$$g(x) =
\int_{-\iy}^{\iy}\frac{1}{\sqrt{8}\pi}\int_{-\iy}^{\iy}
\frac{\hat{f}(\xi) }{\sqrt{i\xi + 0}}
e^{iu^{2}/4\xi}e^{-ixu}d\xi du.$$
Now
\[\int_{-\iy}^{\iy}e^{i({u^2\ov 4\xi}-xu)}\,du
=\sqrt{4\pi}\sqrt{i\xi+0}\,e^{-ix^2\xi},\]
and therefore our symbol is given by
\[ g(x) = \int_{-\iy}^{\iy}\hat{f}(\xi)
e^{-ix^2\xi}d\xi = f(-x^{2}).\]
Thus $P K_{f} P=W(g)$, and the connection to
Wiener-Hopf operators is now apparent.\sp

Lemma 3.1 told us that the operator $U^{-1}(I-P) \A M_{f_{\al}} \A
P U$ converges in the trace norm to $(I-P) W(g) P.$ The next lemma
concerns the strong convergence of the operator
\[B_{\al}(f)=U^{-1}\A M_{f_{\alpha}} \A U.\]
This is the last technical lemma before we can put the pieces together.

\begin{lemma}\label{strong.limit}
The operator $B_{\al}(f)$ converges strongly to $K_{f}$ as $\al\ra\iy$.
\end{lemma}
\proof
We have to show that for any $\ph \in L^2(\R)$
\[\int_{-\iy}^{\iy}\int_{-\iy}^{\iy}{\hat
f(\xi)\ov\sqrt{i\xi+0}}\,e^{-i\xi^3/12\al^3}\,
e^{-i(x+y)\xi/2\al^{3/2}}\,e^{i(x-y)^2/4\xi}\,\ph(y)\,dy\,d\xi\]\[\ra
\int_{-\iy}^{\iy}\int_{-\iy}^{\iy}{\hat
f(\xi)\ov\sqrt{i\xi+0}}\,e^{i(x-y)^2/4\xi}\, \ph(y)\,dy\,d\xi.\]
in $L^2(\R)$. We can restrict ourselves to a dense subset of $\ph$s
since the $B_{\al}(f)$ have uniformly bounded norms.

Write the double integral on the left as \[\int_{-\iy}^{\iy}{\hat
f(\xi)\ov\sqrt{i\xi+0}}\,e^{-i\xi^3/12\al^3}\,e^{-ix\xi/2\al^{3/2}}\,
\int_{-\iy}^{\iy}e^{-iy\xi/2\al^{3/2}}\,e^{i(x-y)^2/4\xi}\,\ph(y)\,dy\,d\xi.\]
This
minus its
purported $L^2$ limit equals the sum of the two error integrals
\bq\int_{-\iy}^{\iy}{\hat
f(\xi)\ov\sqrt{i\xi+0}}\,e^{-i\xi^3/12\al^3}\,e^{-ix\xi/2\al^{3/2}}\,
\int_{-\iy}^{\iy}(e^{-iy\xi/2\al^{3/2}}-1)\,e^{i(x-y)^2/4\xi}\,\ph(y)\,dy
\,d\xi\label{int1}\eq
and \bq\int_{-\iy}^{\iy}{\hat
f(\xi)\ov\sqrt{i\xi+0}}\,(e^{-i\xi^3/12\al^3}\,e^{-ix\xi/2\al^{3/2}}-1)
\,\int_{-\iy}^{\iy}e^{i(x-y)^2/4\xi}\,\ph(y)\,dy \,d\xi.\label{int2}\eq

The operator with kernel $e^{i(x-y)^2/\xi}$ is unitarily equivalent to,
and
therefore has  the same norm as, the operator with kernel
$|\xi|^{1/2}\,e^{i(x-y)^2}$. Thus it has norm $O(|\xi|^{1/2})$. The
function \[(e^{-iy\xi/\al^{3/2}}-1)\,\ph(y)\] has norm
$O(\xi/\al^{3/2})$, assuming
as we
may that $y\,\ph(y)\in L^2$, and it follows that the inner integral in
(\ref{int1}) has norm $O(|\xi|^{3/2}/\al^{3/2})$. Hence (\ref{int1})
itself
has norm
at most $O(1/\al^{3/2})$.

As for (\ref{int2}), the inner integral equals a function
$\psi_{\xi}(x)$ whose
$L^2$ norm is $O(|\xi|^{1/2})$. Write (\ref{int2}) as the sum
\bq\int_{-\iy}^{\iy}{\hat
f(\xi)\ov\sqrt{i\xi+0}}\,(e^{-i\xi^3/12\al^3}-1)\,
e^{-ix\xi/2\al^{3/2}}\,\psi_{\xi}(x)\,d\xi +\int_{-\iy}^{\iy}{\hat
f(\xi)\ov\sqrt{i\xi+0}}\,(e^{-ix\xi/2\al^{3/2}}-1)\,\psi_{\xi}(x)\,d\xi.\label{int3}\eq
The $L^2$ norm of the function
$(e^{-i\xi^3/12\al^3}-1)\,e^{-ix\xi/2\al^{3/2}}\,\psi_{\xi}(x)$ is
$O(|\xi|^{7/2}/\al^3)$
and so the first integral in (\ref{int2}) is a function whose norm is
$O(1/\al^3)$. As for the second integral, observe that
\[\nm(e^{-ix\xi/\al}-1)\,\psi_{\xi}(x)\nm\]
is $O(\xi)$ for all $\xi$ and
$\al$ and tends to 0 as $\al\ra\iy$ for each $\xi$ (by the dominated
convergence
theorem). Hence the integral obtained by taking norm under the integral
sign in
the second integral tends to 0 as $\al\ra\iy$, again by the dominated
convergence theorem. This establishes the claimed strong convergence.
\qed\sp

Now we are ready to begin the final steps in proving (\ref{int.2}). The
operator
$(A_{\alpha}(f))^{n} = (M_{f_{\alpha}}\A P\A)^{n}$ has the same trace
as $(P\A M_{f_{\alpha}}\A P)^{n}$ which in turn has the same trace as
$(P U^{-1}\A M_{f_{\alpha}}\A U P)^{n}=(P B_{\al}(f)P)^{n}$.
In fact for any analytic function $F$
defined on the neighborhood of the spectra of both operators and
satisfying $F(0)=0$ we have
$$ \tr F(A_{\alpha}(f))=\tr F(PB_{\al}(f)P).$$
It is the asymptotics of this last trace we
shall compute. We think of our operators as acting on
$L^{2}(\R^{+}).$\sp

In the following two lemmas $\la$ will be in the resolvent set
of the Wiener-Hopf operator $P K_{f}P=W(g)$. By Lemma~2.3b this implies
that $\la$ is not
in the range of $g(x)=f(-x^2)$, so $f\inv(\{\la\})$
is a compact subset of $(0,\,\iy)$. We can find a Schwartz function
$\tf$ which never
takes the value $\la$ and which equals $f$ outside some larger compact
subset of $(0,\,\iy)$,
and we can find one $\tf$ which serves for all $\la$ in any given
compact subset of
the resolvent set of $W(g)$. This will be our notation in what follows.

\begin{lemma}
Let $\la$ be in the resolvent set of $P K_{f}P=W(g)$.
Then $\lambda$ is also in the resolvent set of $P B_{\al}(f) P$
for sufficiently large $\alpha$ and the inverses have uniformly
bounded norms for $\la$ lying in any given compact subset of the
resolvent set.
\end{lemma}
\proof
Observe that for Schwartz functions $f_1$ and $f_2$, since
$B_{\al}(f_1)\,B_{\al}(f_2)=B_{\al}(f_1\,f_2)$,
\[PB_{\al}(f_1)PB_{\al}(f_2)BP-PB_{\al}(f_1\,f_2)BP=
PB_{\al}(f_1)(P-I)B_{\al}(f_2)BP\]
\bq\ra
PK_{f_{1}}(P-I)K_{f_{2}}P=PK_{f_{1}}PK_{f_{2}}P-PK_{f_{1}f_{2}}P
=W(g_1)\,W(g_2)-W(g_1\,g_2)\label{Blim}\eq
in trace norm since $(I-P)B_{\al}(f_2)BP\ra (I-P)K_{f_{2}}P$ in trace
norm by Lemma \ref{3.1} and
$PB_{\al}(f_1)\ra PK_{f_{1}}$
strongly by Lemma \ref{strong.limit}. This also holds if the $f_i$ are
constants plus Schwartz functions.

We take $f_1=f-\la$ and $f_2=(\tf-\la)\inv$. Observe that
the ``$g$'' corresponding to $f_2$ is $(g-\la)\inv$, so
in this case the relation (\ref{Blim}) reads
\[PB_{\al}(f-\la)PB_{\al}((\tf-\la)\inv)P-PB_{\al}((f-\la)\,(\tf-\la)\inv)P
\ra W(g-\la)\,W((g-\la)\inv)-P.\]
Now
\[PB_{\al}((f-\la)\,(\tf-\la)\inv)P=P+PB_{\al}((f-\la)\,(\tf-\la)\inv-1)P=P+o_1(1),\]
by Lemma 3.2. We conclude that
\bq
PB_{\al}(f-\la)PB_{\al}((\tf-\la)\inv)P=W(g-\la)\,W((g-\la)\inv)+o_1(1).\label{appr}\eq
The analogous formula holds when $f-\la$ and $(\tf-\la)\inv$ are
interchanged. The Wiener-Hopf operators $W(g-\la)$ and $W((g-\la)\inv)$
are
invertible, the first by assumption and the second since
$\ind(g-\la)\inv=-\ind(g-\la)
=0$. The
norms of the inverses are bounded uniformly in $\al$ and $\la$ lying in
a compact set and
the operators $PB_{\al}((\tf-\la)\inv)P$ are uniformly bounded. This
completes the proof.\qed

\begin{lemma}
Suppose $F$ is analytic in a neighborhood of the spectrum of $W(g)$.
Then we have as $\al\ra\iy$
\bq
F(P B_{\al}(f) P) - P B_{\al}(F\circ \tf)P  = F(W(g)) - W(F\circ g)
+o_{1}(1). \label{3.3}
\eq
\end{lemma}
\proof
By Lemma 3.5 $PB_{\al}(f-\la)P$ is
invertible for sufficiently large $\alpha$
with uniformly bounded norm for $\la$ lying in a compact set in the
resolvent set of $W(g).$
Let $\la$ also be in the domain of $F$. Then
\[(PB_{\al}(f-\la)P)\inv-PB_{\al}((\tf-\la)\inv)P=(PB_{\al}(f-\la)P)\inv\,
[I-W(g-\la)\,W((g-\la)\inv)]
+o_1(1)\]\[=W(g-\la)\inv\,[I-W(g-\la)\,W((g-\la)\inv)]+o_1(1).\]
The first equality  follows from (\ref{appr}) and the uniformity of the
norms
of the inverses. The second equality uses the strong convergence of
$(PB_{\al}(f-\la)P)\inv$ to
$W(g-\la)\inv$ and
the fact that $I-W(g-\la)\,W((g-\la)\inv)$ is trace class. Thus
\[(PB_{\al}(f-\la)P)\inv-PB_{\al}((\tf-\la)\inv)P=W(g-\la)\inv-W((g-\la)\inv)+
o_1(1).\]
Multiplying by $F(\la)$ and integrating over an appropriate contour
gives
\[F(PB_{\al}(f)P)-PB_{\al}(F\circ \tf)P=F(W(g))-W(F\circ g)+o_1(1)\]
for any $F$ analytic in a neighborhood of the spectrum of $W(g)$.
\qed\sp

We will be interested in the trace of the first operator on the left
side of (\ref{3.3}).
The next lemma will tell us the trace of the second operator.

\begin{lemma} For any Schwartz function $f$ we have
\[\tr PB_{\al}(f)P=\frac{\al^{3/2}}{\pi}\int_{0}^{\iy}\sqrt{x}f(-x)\,dx
+o(1).\]
\end{lemma}
\proof The kernel of $PB_{\al}(f)P$ equals
$$\frac{1}{\sqrt{8}\pi}\int_{-\iy}^{\iy}\frac{\hat{f}(
\xi)}{\sqrt{i\xi+0}}\,e^{-i\xi^3/12\al^{3}}\,e^{-i(x+y)\xi/2\al^{3/2}}\,
e^{i(x-y)^2/4\xi}\,d\xi\, \ch_{\R^{+}}(x)\ch_{\R^{+}}(y),$$
and thus
$$\tr
PB(f)P=\frac{1}{\sqrt{8}\pi}\int_{0}^{\iy}\int_{-\iy}^{\iy}\frac{\hat{f}(
\xi)}{\sqrt{i\xi+0}}\,e^{-i\xi^3/12\al^{3}}\,e^{-ix\xi/\al^{3/2}}\,
 \,d\xi\,dx .$$
 We write this as

$$\frac{\al^{3/2}}{\sqrt{8}\pi}\int_{0}^{\iy}\int_{-\iy}^{\iy}\frac{\hat{f}(
\xi)}{\sqrt{i\xi+0}}\,e^{-i\xi^3/12\al^{3}}\,e^{-ix\xi}\,
 \,d\xi\,dx ,$$
 and then replace the term $e^{-i\xi^3/12\al^{3}}$ by $1$ just as in
 Lemma \ref{3.1} to find that the trace is given by

$$\frac{\al^{3/2}}{\sqrt{8}\pi}\int_{0}^{\iy}\int_{-\iy}^{\iy}\frac{\hat{f}(
\xi)}{\sqrt{i\xi+0}}\,e^{-ix\xi}\,
 \,d\xi\,dx  + O(\al^{-3/2}).$$
 We can write this in a more familiar form by replacing the term
 $1/{\sqrt{i\xi+0}}$ in the above integral with
 $\frac{1}{\sqrt{\pi}}\int_{-\iy}^{\iy}e^{-iu^{2}\xi}du.$
 Integrating over $\xi$ we find that this equals
 $$\frac{\al^{3/2}}{2\pi}
 \int_{0}^{\iy}\int_{-\iy}^{\iy}f(-u^{2}-x)du\,dx+o(1)$$
 or
 $$\frac{\al^{3/2}}{\pi}\int_{0}^{\iy}\sqrt{x}\,f(-x)\,dx +o(1)$$
 as claimed.\qed\sp

We now derive our main result on determinants of Airy operators which
gives the
 promised formula (\ref{int.2}) for the asymptotics. We assume, as
always, that
 $f$ is a Schwartz function.

 \begin{theorem}
 Assume $g(x)=f(-x^2)\neq-1$. Then as $\al\ra\iy$
 \bq \det\,(I+A_\al(f))= \exp\left\{ c_{1}\,\al^{3/2}+ c_{2}
+o(1)\right\},\label{detA}\eq
where
 $$c_{1} =\frac{1}{\pi}\int_{0}^{\iy}\sqrt{x}\,\log (1+f(-x))\,dx,$$
$$c_{2}=
 \frac{1}{2}\int_{0}^{\iy}x \left( (\log (1+g))^{\check{}}(x)
\right)^{2} dx.$$
 \end{theorem}
\proof Assume first that $\|g\|_{\iy}<1$. Then by Lemma 2.3a the
spectrum of $W(g)$ lies in
the open unit disc with center 0. Therefore $F(z)=\log\,(1+z)$ (the
branch equal to 0 when
$z=0$) is analytic
on the spectrum and so we may apply Lemma 3.6. This and Lemma 3.7 tell
us
 that there is an asymptotic formula of the form (\ref{detA})
 where $c_1$ is as stated and
 \[c_2=\tr [\log\,(I+W(g)) - W(\log\,(1+g))].\]
It is known that this equals the expression given for $c_2$ in the
statement of the
theorem \cite{W2}.

 To remove the restriction on $g$ we introduce a parameter $t$ and we
would like to define
 a family of functions $f_t$ by $1+f_t=e^{t\,\log\,(1+f)}$, so that
 for small enough $t$ our asymptotic formula holds. The problem is that
$-1$ may lie
 in the range if $f$, and even if it didn't we might not have
$i(1+f)=0$,
 which is what we need to
 define a logarithm which is a Schwartz function. So, as in the
preceding
 lemmas, we introduce a function $\tf$ which equals $f$ except on a
compact subset of
 $(0,\,\iy)$ such that $1+f\ne0$ and $\ind\,(1+\tf)=0$. Then we define
$f_t$ by
 \[1+f_t=e^{t\,\log\,(1+\tf)}+f-\tf.\]
 Of course $f_1=f$. Moreover, with $g_t(x)=f_t(-x^2)$,
 \[1+g_t=e^{t\,\log\,(1+g)}\]
 for all $t$. For sufficiently small $t$ we have $\|g_t\|_{\iy}<1$
 so that our formulas hold.

 Observe
 that $\det\,(I+A_\al(f_t))$ is a family of entire functions of $t$
depending on the
 parameter $\al$. Suppose we can show that
 \bq \det\,(I+A_\al(f_t))=O(e^{t\,c_{1}\al^{3/2}})\label{detest}\eq
 for large $\al$ uniformly on compact $t$-sets. Then the limit relation

\[\lim_{\al\ra\iy}e^{-t\,c_1\al^{3/2}}\,\det\,(I+A_{\al}(f_t))=e^{t\,c_2},\]
 which we know holds for sufficiently small $t$, will hold for all $t$
and therefore $t=1$.

 To prove (\ref{detest}) we go back to the $PB_{\al}(f)P$ regarded as
 operators on $L^2(\bf{R}^+)$. We have
 \[\det\,(I+A_\al(f_t))=\det\,PB_{\al}(1+f_t)P.\]
Now $i(1+g_t)=0$ and so $W(1+g_t)$ is invertible by Lemma 2.3b.
Therefore by Lemma~3.5
with $F(z)=z\inv$ we know that $PB_{\al}(1+f_t)P$ will be invertible if
$\al$ is large
enough. (This will hold for all
$t$ in any given compact set.) For these $\al$ we have
\[{d\ov dt}\log\,\det\,(I+A_\al(f_t))=\tr [(PB_{\al}(1+f_t)P)\inv\,
{d\ov dt}
PB_{\al}(1+f_t)P]=\tr[(PB_{\al}(1+f_t)P)\inv\,PB_{\al}(h_t)P],\]
where
\[h_t=\log\,(1+\tf)\,e^{t\,\log\,(1+\tf)}.\]
By Lemma 3.6  with $F(z)=z\inv$ we know that
\[(PB_{\al}(1+f_t)P)\inv=PB_{\al}((1+f_t)\inv)P+O_1(1).\]
Also, by (\ref{Blim}),
\[PB_{\al}((1+f_t)\inv))PB_{\al}(h_t)P=PB_{\al}((1+f_t)\inv\,h_t)P+O_1(1),\]
so that we have shown
\[{d\ov dt}\log\,\det\,(I+A_\al(f_t))=\tr
PB_{\al}((1+f_t)\inv\,h_t)P+O(1).\]
But Lemma 3.7 tells us that with an error $o(1)$
\[\tr PB_{\al}(1+f_t)\inv\,h_t)P=
{\al^{3/2}\ov \pi}\int_0^{\iy}(1+f_t(-x))\inv\,h_t(-x)\,\sqrt{x}\,dx\]
\[={\al^{3/2}\ov \pi}\int_0^{\iy}\,\sqrt{x}\,\log\,(1+f(-x))\, dx.\]
because $\tf=f$ on $\R^-$. Thus,
\[{d\ov dt}\log\,\det\,(I+A_\al(f_t))=
{\al^{3/2}\ov \pi}\int_0^{\iy}\sqrt{x}\,\log\,(1+f(-x))\, dx+O(1).\]
Integrating over $t$ from 0 to $t$ and exponentiating gives
(\ref{detest}) and completes
the proof.\qed

 \section{Applications to random matrices}

 Theorem 3.7
 can be applied to find limiting distribution functions for a class of
random
 variables which are functions of the eigenvalues of a random matrix.
 In many different ensembles of
 matrices it has been shown that the  distribution functions are
 asymptotically normal \cite{F,B,BT,Be}, and this will be shown also to
be the case in
 our examples. The term {\em ensemble} refers to the
 probability density assigned to some space of matrices, and this
 in turn induces a density on the space of eigenvalues of the
 matrices. For the Gaussian Unitary Ensemble (GUE) the density on the
space of
 eigenvalues is given by

 \begin{equation}
  P_N(x_1,\ldots,x_N) = \frac{1}{N!}\det K(x_i,x_j) \left.
\right|_{i,j=1}^N
\end{equation}
where
\begin{equation}
  K_N(x,y)= \sum_{i=0}^{N-1}\phi_i(x) \phi_i(y).
\end{equation}
and $\phi_{i}$ is obtained by orthonormalizing the sequence
$\{x^{i}e^{-x^2/2}\}$
over $\bf{R}.$
If $N$ is large it is also well known that the density of the
eigenvalues
is supported on approximately the interval $(-\sqrt{2N}, \sqrt{2N}).$
These facts can be found in \cite{M}.

The random variables of interest here
are ones that are often called  {\it linear statistics} and are of the
form  $$\sum_{i=1}^{N}f(\la_{i}/\al),$$ where $\la_i$ are the
eigenvalues and $f$
is an appropriate
function. Our goal is to study these random variables applied to the
eigenvalues near the edge of the spectrum and to this end we rescale and
replace
the sum by
$$\sum_{i=1}^{N}f(2^{1/2}N^{1/6}(\la_{i} - \sqrt{2N})/\al).$$
The purpose of the translation by the term $\sqrt{2N}$ is to move to the
edge and
the factor $2^{1/2}N^{1/6}$ has the effect of making the eigenvalue
density of the order 1. Otherwise the eigenvalues
``bunch up'' or ``spread out'' and all the results become more or less
trivial.

To study the distribution function of this random variable we use its
characteristic
function, or
inverse Fourier transform. This characteristic function is given by
\[ \phi_N(s)=
  \int_{-\infty}^{\infty}\cdots \int_{-\infty}^{\infty}
  e^{is\sum_{j=1}^{N}f(2^{1/2}N^{1/6}(x_{j} - \sqrt{2N})/\al))}\,
  P_N(x_1,\ldots,x_N) \,dx_1 \cdots dx_N .
   \]
 It is a general fact that
 \[\int_{-\infty}^{\infty}\cdots \int_{-\infty}^{\infty}
  \prod_{i=1}^N (1+g(x_i))\,
  P_N(x_1,\ldots,x_N)\, dx_1 \cdots dx_N =\det(I+gK_N),\]
  where $g$ denotes multiplication by $g(x)$ and $K_N$ denotes the
operator with kernel
  $K_N(x,y)$. This can be obtained by expanding out the product in the
integrand, using
  the formula
  \[
  \frac{N!}{(N-n)!}\int \cdots \int P_N(x_1,\ldots,x_n,x_{n+1},
  \ldots, x_N) \,dx_{n+1} \cdots dx_N = \det K_{N}(x_i,x_j) \left.
  \right|_{i,j=1}^n\]
for the $n$-point correlation function, and then recognizing the
resulting sum of multiple
integrals as the expansion of the Fredholm determinant. Or it can be
obtained
by a simpler algebraic device \cite{TW1}.
In our case $1+g(x) = \exp\{f(2^{1/2}N^{1/6}(x - \sqrt{2N})/\al)\}$. If
we make the
changes of variable
\[x\ra \frac{x}{2^{1/2}N^{1/6}}+\sqrt{2N},\ \ \ y\ra
\frac{y}{2^{1/2}N^{1/6}}+\sqrt{2N}\]
we find that the characteristic function equals the determinant of $I$
plus the
operator with kernel
\[
(e^{isf(x/\al)}-1)\,\frac{1}{2^{1/2}N^{1/6}}\,K_N(\frac{x}{2^{1/2}N^{1/6}}+
  \sqrt{2N} ,\,\frac { y}{2^{1/2}N^{1/6}} + \sqrt{2N}).\]

Now one has the scaling limit
\[ \lim_{N\ra \iy}
\frac{1}{2^{1/2}N^{1/6}}\,K_N(\frac{x}{2^{1/2}N^{1/6}}+
  \sqrt{2N} ,\,\frac { y}{2^{1/2}N^{1/6}} + \sqrt{2N})\]
\[=  \frac{A(x)A'(y) - A'(x)A(y)}{x-y},\]
precisely the Airy kernel.
Thus we see that the large $N$ limit of the characteristic function
equals
$\phi(s)=\det\,(I + A_{\alpha}(h))$ where $h(x) = e^{isf(x)}-1.$
Our asymptotic formula yields
\[ \phi(s) = \exp{ \left\{
\frac{is\al^{3/2}}{\pi}\int_{0}^{\iy} \,\sqrt{x} f(-x) \,dx -
 \frac{s^{2}}{2}\int_{0}^{\iy}x\, (\check{g}(x))^{2}\,dx
+o(1)\right\}},\]
 where as before $g(x) = f(-x^{2})$.

 Notice that the limiting characteristic
 function is quadratic in $s$ and hence the distribution is
 asymptotically normal. Of course this is not surprising since this
occurs for other matrix
 ensembles and other scaling limits. Notice, though, that in this case
the mean and
 variance of the limiting distribution only depend on the
 negative values of the argument of the original $f.$ This is a
reflection of the fact
 that the Airy function goes rapidly to zero for
 positive values and oscillates and tends to zero slowly for negative
 values of the argument. A question left to the future is
 how the asymptotics of functions $A(x)$ and $B(x)$ in a
 kernel of the form $${A(x)B(y)-A(y)B(x)\ov x-y}$$ affect the
 asymptotics of the corresponding distribution functions.

\end{document}